\documentclass[reqno]{amsart}
\usepackage{amssymb}
\usepackage{epsf}

\newtheorem{theorem}{Theorem}
\newtheorem{proposition}[theorem]{Proposition}

\theoremstyle{definition}

\theoremstyle{remark}
\newtheorem*{remark}{Remark}

\numberwithin{equation}{section}

\begin{document}

\title[Zeros of the Associated Legendre Functions]
{Remarks on the Zeros of the Associated Legendre Functions with
Integral Degree}

\author{J.F. van Diejen}
\address{
Instituto de Matem\'atica y F\'{\i}sica, Universidad de Talca,
Casilla 747, Talca, Chile}

\thanks{Work supported in part by the `Anillo Ecuaciones Asociadas a
Reticulados', financed by the World Bank through the `Programa
Bicentenario de Ciencia y Tecnolog\'{\i}a'.}

\date{February 2007}
\

\begin{abstract}
We present some formulas for the computation of the zeros of the
integral-degree associated Legendre functions with respect to the
order.
\end{abstract}

\maketitle

\section{Introduction}\label{sec1}
The {\em associated Legendre functions} (or {\em spherical
functions}) are given explicitly by
\cite{abr-ste:handbook,grad-ryz:table}
\begin{eqnarray}
\Gamma (1-z) P^z_n(\tanh(x)) &=&
\Bigl(\frac{1+\tanh(x)}{1-\tanh(x)}\Bigr)^{z/2} F
(-n,n+1;1-z; {\textstyle\frac{1-\tanh(x)}{2}}) ,\nonumber \\
&=& \frac{\exp(zx)}{(1+e^{-2x})^n} F(-n,-n-z;1-z;-e^{-2x}),
\label{alf}
\end{eqnarray}
where $\Gamma(-)$ refers to the Euler gamma function and
$F(-,-;-;-)$ to the Gauss hypergeometric series, and where we have
fixed the normalization such that no gamma factors appear in front
of the hypergeometric series representation.

When the {\em degree} $n$ is integral the hypergeometric series on
the first line of Eq. \eqref{alf} terminates, whence it is then
polynomial in the argument $\tanh (x)$ and---with the present
normalization---rational in the {\em order} $z$. The purpose of the
present note is to characterize the locations of zeros of the
integral-degree associated Legendre functions with respect to the
order $z$. In Refs. \cite{dri-dur:zeros,dri-dur:trajectories}
detailed information concerning the locations of the zeros of the
associated Legendre functions with respect to the argument was
provided for the situation that the order and the degree differ by a
positive integer. Previously, a numerical study of the locations of
the zeros with respect to the degree (for integral order and various
values of the argument) was presented in Ref. \cite{bau:tables}.

Questions about the zeros of the associated Legendre functions fit
within a rich tradition of research concerning the locations of the
zeros of {\em orthogonal} polynomials
\cite{sze:orthogonal,ism:classical}. An important difference is,
however, that in the present context we are not in the position to
exploit an orthogonality structure. Instead, the idea behind the
methods below is to employ a connection with the (inverse)
scattering theory of the one-dimensional Schr\"odinger equation with
P\"oschl-Teller potential \cite{flu:practical}, which reveals that
the integral-degree associated Legendre functions correspond to
reflectionless wave functions that can be expressed in terms of the
tau functions of the Korteweg-de Vries hierarchy
\cite{die-kir:combinatorial}. Detailed information on the zeros can
then be obtained from Refs. \cite{die:zeros,die-pus:reflectionless},
where the behavior of the zeros of such reflectionless wave
functions was studied with the aid of (Ruijsenaars-Schneider type)
integrable particle systems (cf. also Ref.
\cite{erc-lev-zha:behavior}, where information on the zeros of
reflectionless Schr\"odinger wave functions was obtained via inverse
scattering techniques).

This note is organized as follows. First we present a system of
algebraic equations for the zeros of the associated Legendre
function in Section \ref{sec2}. Next, it is shown in Section
\ref{sec3} that these algebraic equations imply a nonlinear system
of coupled differential equations. The solution of these
differential equations formulated in Section \ref{sec4} then
provides us with explicit information on the locations of the zeros
of the integral-degree associated Legendre functions.

\section{Algebraic Equations}\label{sec2}
Let us abbreviate the renormalized associated Legendre function in
Eq. \eqref{alf} by $\psi_n(x,z)$. Upon writing the hypergeometric
series explicitly:
\begin{equation}\label{alfn}
\psi_n(x,z) =\frac{\exp(zx)}{(1+e^{-2x})^n} \sum_{m=0}^n
e^{-2mx}\binom{n}{m} \prod_{j=1}^m \frac{z+n+1-j}{z-j} ,
\end{equation}
it is seen that the function in question admits a factorization of
the form
\begin{equation}\label{factor}
\psi_n(x,z)= \frac{\exp(zx)}{(1+e^{-2x})^n} \prod_{j=1}^n
\frac{z-z_j(x)}{z-j} .
\end{equation}

The following proposition provides a (Bethe type) system of
algebraic equations for the zeros $z_1(x),\ldots ,z_n(x)$.
\begin{proposition}[Algebraic System]\label{bethe:prp}
The zeros $z_1(x),\ldots ,z_n(x)$ of the normalized integral-degree
associated Legendre function $\psi_n(x,z)=\Gamma (1-z)
P^z_n(\tanh(x))$ satisfy the rational system
\begin{equation}\label{bethe:eq}
\prod_{j =1}^n \frac{\ell-z_j (x)}{\ell+z_j (x)} = (-1)^{n-\ell}
\exp(-2\ell x),\qquad \ell =1,\ldots ,n.
\end{equation}
\end{proposition}
\begin{proof}
It follows from the Gauss hypergeometric equation (and can also be
inferred directly with the aid of Eq. \eqref{alfn}) that the
associated Legendre function $\psi_n(x,z)$ solves the Schr\"odinger
equation on the line with P\"oschl-Teller potential
\begin{equation*}
 \Bigl(
\frac{\text{d}^2}{\text{d}x^2}+\frac{n(n+1)}{\cosh^2(x)}-z^2\Bigr)
\psi_n(x,z)=0.
\end{equation*}
From Eq. \eqref{alfn} it is manifest that $\psi_n(x,z)\to \exp (zx)$
for $x\to+\infty$ and $\psi_n(x,z)\to \exp (zx)\prod_{j=1}^n
\frac{z+j}{z-j}$ for $x\to -\infty$. In other words, $\psi_n(x,z)$
provides a Jost solution of the Schr\"odinger eigenvalue problem in
question and the pertinent P\"oschl-Teller potential
$n(n+1)/\cosh^2(x)$ is {\em reflectionless}. We read-off from the
asymptotics for $x\to\pm\infty$ that the bound states correspond to
the spectral values $z=-j$, $j=1,\ldots,n$; the associated
normalization constants are given explicitly by (cf.
\cite[\text{Eqs}. (7.122.1), (8.737.1,2)]{grad-ryz:table})
\begin{equation}\label{norm}
\nu_j=\left( \int_{-\infty}^\infty
\psi^2_n(x,-j)\text{d}x\right)^{-1}= j\binom{2j}{j}\binom{n+j}{n-j}
= (-1)^{n-j}2j\prod_{k=1,k\neq j}^n \frac{j+k}{j-k}.
\end{equation}
The proposition now follows upon specialization of a system of Bethe
type equations for the zeros of reflectionless Jost functions
derived in Refs. \cite{die:zeros,die-pus:reflectionless}:
\begin{equation*}
2\kappa_\ell \prod_{j=1}^n
\frac{\kappa_\ell-z_j(x)}{\kappa_\ell+z_j(x)} \prod_{j=1,j\neq
\ell}^n \frac{\kappa_\ell+\kappa_j}{\kappa_\ell-\kappa_j}=\nu_\ell
\exp(-2x\kappa_\ell),\qquad \ell=1,\ldots ,n,
\end{equation*}
where $-\kappa_1,\ldots ,-\kappa_n$ represesent the spectral values
parametrizing the discrete spectrum and $\nu_1,\ldots ,\nu_n$ denote
the values of the associated normalization constants (so in our case
$\kappa_j=j$ and $\nu_j$ is given by the expression in Eq.
\eqref{norm}).
\end{proof}

\section{Differential Equations}\label{sec3}
By  differentiating the algebraic equations in Proposition
\ref{bethe:prp}, one arrives at a (Dubrovin type) system of
first-order differential equations for the zeros.

\begin{proposition}[First-order System]\label{dubrovin:prp}
The zeros $z_1(x),\ldots ,z_n(x)$ of the normalized integral-degree
associated Legendre function $\psi_n(x,z)=\Gamma (1-z)
P^z_n(\tanh(x))$ satisfy the system of first-order differential
equations
\begin{subequations}
\begin{equation}
z_\ell^\prime (x) = \frac{\prod_{j =1}^n (j^2-z_\ell^2(x))}{\prod_{j
=1, j\neq \ell}^n (z_j^2(x)-z_\ell^2(x))} ,\qquad \ell=1,\ldots ,n,
\end{equation}
with the initial condition
\begin{equation}
 \{ \lim_{x\to 0}z_1(x),\ldots ,\lim_{x\to 0}z_n(x)\} = \{ n-1,n-3,\ldots, -(n-3),
-(n-1) \} .
\end{equation}
\end{subequations}
\end{proposition}
\begin{proof}
When taking the (logarithmic) derivative of the algebraic equations
in Proposition \ref{bethe:prp}, one arrives at the identities
\begin{equation*}
\sum_{j=1}^n \frac{z_j^\prime (x)}{\ell^2-z_j^2(x)}=1,\qquad
\ell=1,\ldots ,n.
\end{equation*}
These identities  may be thought of as a system of $n$ linear
equations for the derivatives $z_1^\prime (x),\ldots ,z_n^\prime
(x)$. Solving this linear system leads to the differential equations
stated in the proposition. For $x=0$, the hypergeometric series in
the second line of Eq. \eqref{alf} simplifies:
\begin{equation*}
\psi_n(0,z)=\Gamma (1-z)P_n^z(0)=2^{-n} F(-n,-n-z;1-z;-1)=
\prod_{j=1}^n \frac{z+n+1-2j}{z-j},
\end{equation*}
whence for $x\to 0$ the zeros of the integral-degree associated
Legendre function tend to $n+1-2j$, $j=1,\ldots ,n$.
\end{proof}

\begin{remark}
Upon substituting the product representation of the form in Eq.
\eqref{factor} into the Schr\"odinger equation with P\"oschl-Teller
potential in the proof of Proposition \ref{bethe:prp}, it is seen
that the zeros of the normalized integral-degree associated Legendre
function $\psi_n(x,z)=\Gamma (1-z) P^z_n(\tanh(x))$ also satisfy the
following system of (Ruijsenaars-Schneider type) second-order
differential equations (cf. Refs.
\cite{die:zeros,die-pus:reflectionless})
\begin{equation}
z_\ell^{\prime\prime}(x)+2z_\ell(x)z_\ell^\prime(x) = \sum_{j =1,
j\neq \ell}^n \frac{2z_\ell^\prime (x)z_j^\prime(x)}{z_\ell(x)-z_j
(x)},\qquad \ell=1,\ldots ,n.
\end{equation}
From Proposition \ref{dubrovin:prp} it follows that the
corresponding initial conditions are given by
\begin{subequations}
\begin{eqnarray}
z_\ell (0)&=&n+1-2\ell, \\
z_\ell^\prime (0) &=&\frac{\displaystyle -
\prod_{\begin{subarray}{c}j=1\\j\neq |n+1-2\ell|\end{subarray}}^n
((n+1-2\ell)^2-j^2)}{\displaystyle 2^{n-1}\prod_{\begin{subarray}{c}j=1\\
j\neq \ell,n+1-\ell\end{subarray}}^n 2(\ell-j)^2} ,
\end{eqnarray}
$\ell=1,\ldots ,n$ (where we have ordered the zeros from large to
small).
\end{subequations}
\end{remark}

\section{Trajectories}\label{sec4}
We conclude by providing a characterization of the zeros of the
associated Legendre functions in terms of the eigenvalues of an
explicit square matrix $\mathbf{Z}$ of size equal to the degree $n$.

\begin{proposition}[Trajectories]\label{trajectories:prp}
The zeros $z_1(x),\ldots ,z_n(x)$ of the normalized integral-degree
associated Legendre function $\psi_n(x,z)=\Gamma (1-z)
P^z_n(\tanh(x))$ are given by the (real) eigenvalues of the matrix
\begin{subequations}
\begin{equation}
\mathbf{Z}(x):=\mathbf{K}(\mathbf{I}-\exp(-2x\mathbf{K})\mathbf{N})
(\mathbf{I}+\exp(-2x\mathbf{K})\mathbf{N})^{-1},
\end{equation}
where $\mathbf{K}:=\text{diag}(1,2,\ldots ,n)$, $\mathbf{N}$ is the
$n\times n$ matrix with components
\begin{equation}
\mathbf{N}_{j,k}:= \binom{2j}{j}\binom{n+j}{n-j}\frac{j}{j+k},\qquad
1\leq j,k\leq n
\end{equation}
\end{subequations}
(and $\mathbf{I}$ represents the $n$-dimensional identity matrix).
\end{proposition}
\begin{proof}
Let $\tilde{\mathbf{Z}}:=\text{diag}(z_1(x),\ldots ,z_n(x))$ and let
$\mathbf{U}$ denote the orthogonal matrix with components
\begin{equation*}
\mathbf{U}_{j,k}:= \left(\frac{\prod_{\ell =1}^n
(\ell^2-z_j^2(x))}{\prod_{\ell =1,\ell\neq j}^n
(z_\ell^2(x)-z_j^2(x))}\right)^{\frac{1}{2}} \frac{1}{z_j^2(x)-k^2}
\left(\frac{\prod_{\ell =1}^n (k^2-z_\ell^2(x))}{\prod_{\ell
=1,\ell\neq k}^n (k^2-\ell^2)}\right)^{\frac{1}{2}},
\end{equation*}
$1\leq j,k\leq n$. Then the matrix
\begin{equation*}
\tilde{\mathbf{N}}:=
(\mathbf{I}-\mathbf{K}^{-1/2}\mathbf{U}^{-1}\tilde{\mathbf{Z}}\mathbf{U}\mathbf{K}^{-1/2})
(\mathbf{I}+\mathbf{K}^{-1/2}\mathbf{U}^{-1}\tilde{\mathbf{Z}}\mathbf{U}\mathbf{K}^{-1/2})^{-1}
\end{equation*}
has components
\begin{equation*}
\tilde{\mathbf{N}}_{j,k}
=\frac{\tilde{\nu}_j^{1/2}\tilde{\nu}_k^{1/2}}{j+k},\qquad 1\leq
j,k\leq n,
\end{equation*}
where $\tilde{\nu}_\ell:=2\ell \prod_{j=1}^n
\frac{\ell-z_j(x)}{\ell+z_j(x)} \prod_{j=1,j\neq \ell}^n
\frac{\ell+j}{\ell-j}\stackrel{\text{Eq.
\eqref{bethe:eq}}}{=}\ell\binom{2\ell}{\ell}\binom{n+\ell}{n-\ell}\exp(-2\ell
x)$, $\ell =1,\ldots ,n$. By inverting the relation between
$\tilde{\mathbf{Z}}$ and $\tilde{\mathbf{N}}$ it is readily seen
that the zeros $z_1(x),\ldots ,z_n(x)$ are given by the eigenvalues
of the symmetric matrix
\begin{equation*}
\mathbf{K}^{1/2}
(\mathbf{I}-\tilde{\mathbf{N}})(\mathbf{I}+\tilde{\mathbf{N}})^{-1}\mathbf{K}^{1/2},
\end{equation*}
which differs from the matrix $\mathbf{Z}(x)$ formulated in the
proposition by a diagonal similarity transformation.
\end{proof}

\begin{figure}[t]
\centerline{\epsffile{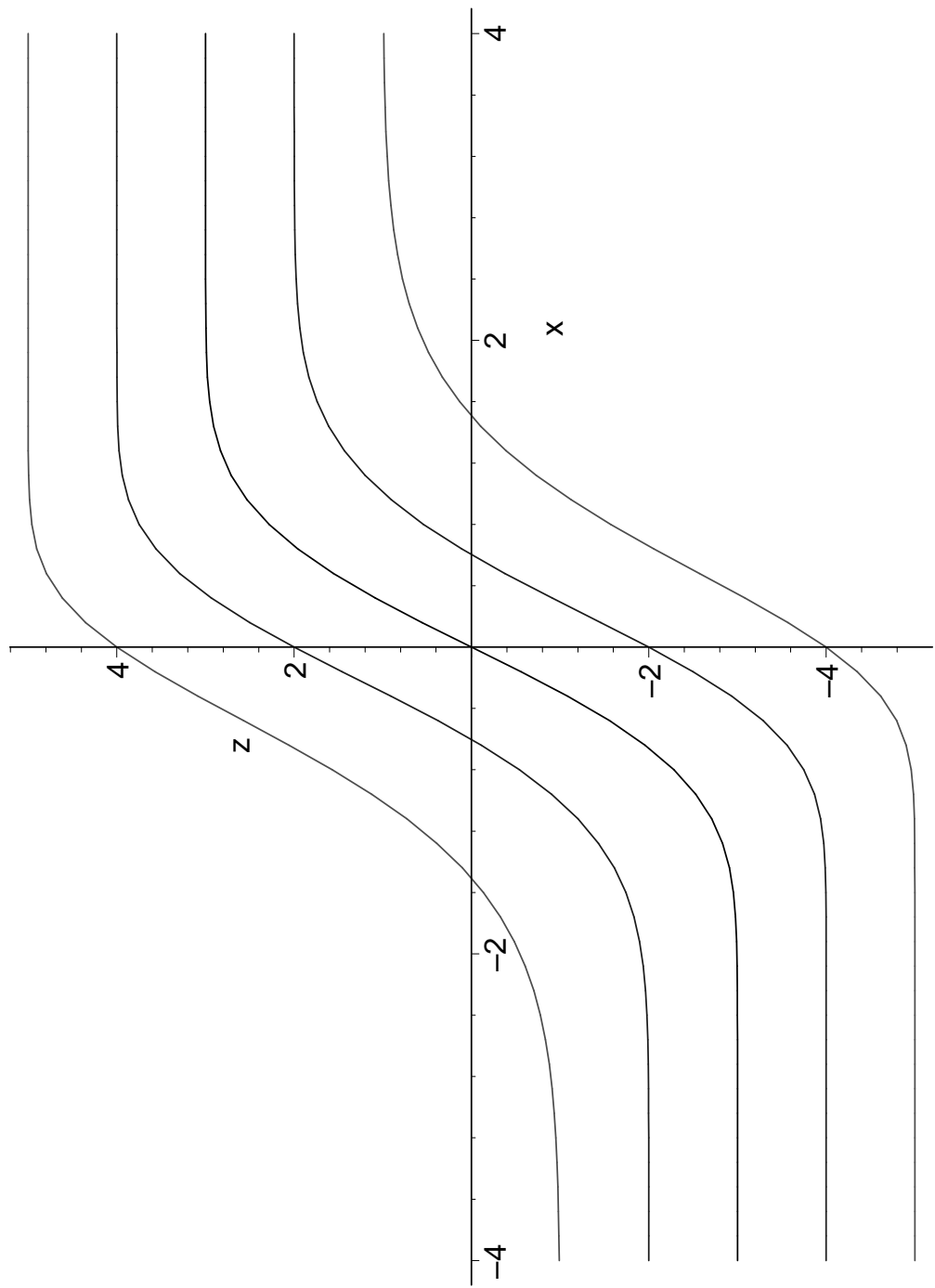}} \caption{Trajectories of the zeros
of $\Gamma (1-z)P_n^z(\tanh(x))$ with respect to $z$ as a function
of $x\in\mathbb{R}$ for $n=5$.} \label{fig1}
\end{figure}

\begin{figure}[t]
\centerline{\epsffile{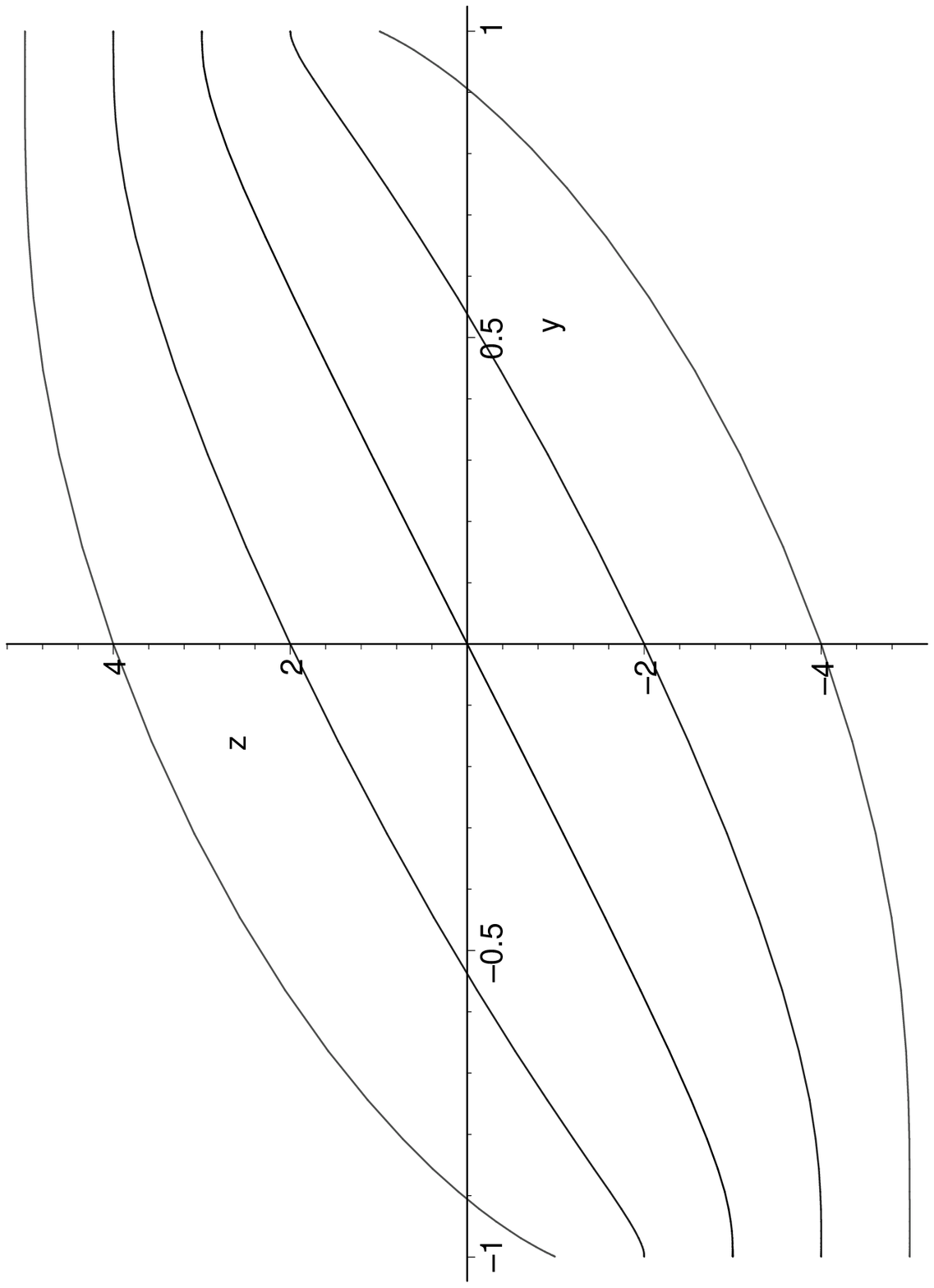}} \caption{Trajectories of the zeros
of $\Gamma (1-z)P_n^z(y)$ with respect to $z$ as a function of
$-1\leq y\leq 1$ for $n=5$. } \label{fig2}
\end{figure}

When $x$ runs from $-\infty$ to $+\infty$ along the real axis the
argument $y=\tanh(x)$ of the associated Legendre function $P_n^z(y)$
varies from $-1$ to $1$. It is clear from Proposition
\ref{trajectories:prp} that its zeros $\{ z_1(x),\ldots ,z_n(x)\}$
move in this situation from $\{ -1,-2,\ldots ,-n\}$ to $\{
1,2,\ldots ,n\}$. It moreover follows from the general analyis in
Ref. \cite{die-pus:reflectionless} that the corresponding
trajectories of the zeros are analytic in $x$ and strictly
monotonously increasing in such a way that crossings do not occur
(i.e. the zeros remain simple for all $x\in\mathbb{R}$). To
illustrate this state of affairs, in Figures \ref{fig1} and
\ref{fig2} the trajectories of the zeros are plotted for $n=5$ as
function of $x$ and of $y=\tanh(x)$, respectively.

\begin{remark}
It is manifest from Eq. \eqref{factor} and Proposition
\ref{trajectories:prp} that the integral-degree associated Legendre
function admits the following determinantal representation
\begin{subequations}
\begin{equation}
\Gamma (1-z) P^z_n(\tanh(x))=\frac{e^{zx}\det
(z\mathbf{I}-\mathbf{Z}(x))}{(1+e^{-2x})^n\prod_{j=1}^n(z-j)},
\end{equation}
with the matrix  $\mathbf{Z}(x)$ as defined in Proposition
\ref{trajectories:prp}, or equivalently:
\begin{eqnarray}
\det (z\mathbf{I}-\mathbf{Z}(x)) &=&
F(-n,-n-z;1-z;-e^{-2x})\prod_{j=1}^n(z-j) \nonumber \\
&=& \sum_{m=0}^n e^{-2mx}\binom{n}{m} \prod_{j=1}^m (z+n+1-j)
\prod_{k=m+1}^n(z-k) .
\end{eqnarray}
\end{subequations}
\end{remark}

\bibliographystyle{amsplain}

\end{document}